\documentclass[12pt,a4paper]{article}
\usepackage{t1enc}
\usepackage[dvips,final]{graphicx} 
\usepackage{tabularx} 
\usepackage{amsfonts}
\usepackage{amssymb,amsmath}

\def\du{\unskip\smash{\lower 1.4ex \hbox{\char34}}\kern-.2ex}
\def\hu{\kern-.2ex\hbox{\char92}}

\def\XXint#1#2#3{{\setbox0=\hbox{$#1{#2#3}{\int}$}
     \vcenter{\hbox{$#2#3$}}\kern-.5\wd0}}

\newcommand{\bdis}{\begin{displaymath}}
\newcommand{\edis}{\end{displaymath}}
\newcommand{\be}{\begin{equation}}
\newcommand{\ee}{\end{equation}}
\newcommand{\mbb}{\mathbb}
\newcommand{\mcal}{\mathcal}

\newcommand{\noi}{\noindent}

\newtheorem{pr*}[thm]{*}

\begin{document}
\baselineskip=6mm
\newpage

\title{On computation of clustering coefficient in a class of random networks}
\author{Michal Demetrian\footnote{Comenius University,
Mlynsk\' a Dolina, 842 48, Bratislav IV, Slovakia, e-mail: demetrian@fmph.uniba.sk},
Martin Neh\' ez\footnote{VSM School of Management, City University,
Bratislava, Slovak Republic and FEE CTU, Karlovo N\' am\v est\' i 13, 121 35 Prague, Czech Republic, e-mail:\ nehez841@gmail.com}}
\date{}

\maketitle

\abstract{The random networks enriched with additional structures as metric and group-symmetry in background metric space are investigated.
The important quantities like
he clustering coefficient as well as the mean degree of separation in such networks are effectively computed with help of additional structures.
Representative models are discussed in details.}

\noi
{\bf keywords}:\ random graph, clustering, degree of separation \\
{\bf MCS}:\ 05Cxx

\section{Introduction}

It is well-known for about 40 years that the mean degree of separation of humans (the mean number of people needed to
bind two randomly chosen persons via a chain of acquaintance) is surprisingly small, it takes
the value about six. This was discovered for the first time by the American psychologists Stanley Milgram
\cite{milgram}. Since his
days this very interesting phenomenon attracts the interest of scientists from various fields of research. As shown in \cite{watts}, the
small world phenomenon appears as a generic feature of many natural as well as artificial random graphs. To explore these
properties of the random networks and to understand them one finds often helpful to use methods
developed for statistical mechanics \cite{reka}. Generally, there is no hope to obtain exact results for main quantities describing
random networks of real interest. However, under some additional condition one can build up simplified models obeying some kind of symmetries
that allow for expressing crucial quantities (their approximations) in closed form. Methods of computation can be divided into two classes; one of
them uses purely discrete point of view and the second one, which is much more effective if applicable, uses continuous limit of the network. We shall
be oriented in the latest. \\
The importance of the use of a continuum limit in various discrete random structures is
surely undisputable because the exact discrete mathematics based methods have commonly bounded area of use. In current research the
continuum limit methods are used widely, for the use in direct random graph theory see e.g. \cite{medo}, an interesting result is achieved by
these methods studying some hidden variables (variables that are not measured directly like degree of a node) in a
kind of biological networks \cite{miller}. One of the most important areas of research where these continuum limits are used is the problem
of phase transition in random structures (critical phenomena), for the review see e.g. \cite{lambiotte} and further references therein. \\
The main aim of the present paper is to show how to compute the relevant quantities like the clustering coefficient and the mean
degree of separation of nodes in a given random network that is build up on a metric space or a space where we can define some
oriented distance. The crucial assumption is that the network obeys some symmetry.
The metric serves us to
define the probability that two points (nodes, persons) are linked via a kind of acquaintance, and this probability
is a given function of the distance between the nodes.

\section{Definitions and basic properties of relevant quantities}

Let us consider the number $N$ of entities (people) represented as points lying in the metric space $\Omega$ equipped
with the distance $d$.
The acquaintance (defined in a suitable way) between two people represented by the points $A$ and $B$
is recorded as a link between these points. With the help of the metric considered on $\Omega$ we can construct
a probability $Q(A,B)$ that $A$ knows $B$ (and vice-versa, we consider always the reflexivity of knowing someone) as
a function of the distance between $A$ and $B$:
\be \label{prob}
Q(A,B)=Q(d(A,B)) .
\ee
The simplest quantity characterizing our network that is given by $Q$ is the average number of
acquaintance for any person $A$:
\be \label{nacq}
\mcal{N}=\sum_B Q(A,B) ,
\ee
where the summation runs over all points of considered network. If one assumes that the function $Q$ varies slowly
on the scale of elementary distance in the network, then the above written sum can be approximated by the integral:
\be \label{nacq1}
\mcal{N}=\int_{\Omega} Q(x,y){\rm d}y
\ee
with proper measure ${\rm d}y$ on $\Omega$. (For example, $\Omega$ can be considered to be
an submanifold in some Euclidean space $\mbb{R}^n$ and then it is equipped naturally with both metric and volume form.) \\

Let $b$ be the distance
between the two considered points. Then the average degree of separation can be written as the mean value of the
function (mean with respect to its argument $S$) $P(S,b)$ that is the probability that the degree of separation of
$A$ and $B$ at relative distance $b$ equals $S$. \\

First of all, $P(0,b)$ is given trivially by definition of $Q$ itself:
\be \label{p0}
P(0,b)=Q(b) .
\ee
Other $P(S,b)$'s are to be computed nontrivially. And we are able to obtain analytical results only if some
approximations in computing $P(S,b)$ are applicable. Let us discuss $P(1,b)$. This means that there is just one
person, say $C$ separating the two chosen persons $A$ and $B$. (Let us note that $C$ is necessarily different from $A$
as well as from $B$, otherwise $A$ knows $B$.) With respect what has just been said we obtain the probability in question:
\be \label{p1}
P(1,b)=\sum_{C}Q(A,C)Q(C,B)(1-Q(A,B)) .
\ee
Analogically we comes at the expression for the probability that the separation of the two points at fixed distance is
two:
\begin{eqnarray} \label{p2}
& &
P(2,b)= \sum_{C_1}\sum_{C_2}Q(A,C_1)Q(C_1,C_2)Q(C_2,B)(1-Q(A,C_2))\times \nonumber \\
& &
(1-Q(B,C_1))(1-Q(A,B)) ,
\end{eqnarray}
and we could continue to further separation indices in the same way. However, instead of doing this, we are going to
describe an approximation we shall use in the future considerations. We assume that the probabilities $Q(A,B)$ are
small enough and therefore the product of two such probabilities (that is of the higher order of smallness) can be
neglected in the expression (\ref{p2}) as well as in others $P(S,b)$'s. In this way one obtains instead of (\ref{p2}) its
reduction:
\be \label{p2app}
P(2,b)\approx \sum_{C_1}\sum_{C_2}Q(A,C_1)Q(C_1,C_2)Q(C_2,B) .
\ee
By means of derivation (\ref{nacq1}) from (\ref{nacq}) we can write down for $P(2,b)$ also the following integral
formula
\begin{eqnarray} \label{p2appint}
& &
P(2,b)\approx \nonumber \\
& &
\int_{\Omega\times\Omega}Q(A,x)Q(x,y)Q(y,B)(1-Q(A,y))(1-Q(x,B))(1-Q(A,B)){\rm d}x{\rm d}y \nonumber \\
& &
\approx \int_{\Omega\times\Omega}Q(A,x)Q(x,y)Q(y,B){\rm d}x{\rm d}y .
\end{eqnarray}
Obviously, the same can be done for all other $P(S,b)$'s with the similar result. The special structure of the
integral formula (\ref{p2appint}) will allow for interesting expression for the essential quantity we
shall turn now our attention to.
The quantity that is effectively and often used in analysis of random networks is the mean
clustering coefficient $\langle C\rangle$, \cite{watts2}. The quantity $\langle C\rangle$ is the probability that (any) two
acquaintances of given person $A$ knows each other. By this definition, $\langle C\rangle$ can be computed modifying
equation (\ref{p2appint}) - namely we must identify $A$ and $B$ and avoid interchanging the positions of persons
$C_1$ and $C_2$ that means to reduce the average number of triples by the factor $1/2$:
\bdis
\frac{1}{2}\int_{\Omega\times\Omega}Q(A,x)Q(x,y)Q(y,A){\rm d}x{\rm d}y .
\edis
To obtain $\langle C\rangle$ we have to divide the above written average number of triples by the number of possible
triples that is given by the fact that average number of acquaintances for any person is $\mcal{N}$, so the
normalizing constant we are looking for reads:
$
\frac{\mcal{N}(\mcal{N}-1)}{2}\approx \frac{\mcal{N}^2}{2} ,
$
where we have assumed $\mcal{N}\gg 1$. Finally, the clustering coefficient is given by the following simple equation:
\be \label{ccoef}
\langle C\rangle =\frac{1}{\mcal{N}^2}\int_{\Omega\times\Omega}Q(A,x)Q(x,y)Q(y,A){\rm d}x{\rm d}y .
\ee
The idea of the following is that the above written quantities can be expressed in approximate closed form using additional assumptions on
symmetry of the metric space in use. We shall demonstrate this in details on a concrete example, however, the idea can be extended to wide
range of symmetric spaces in the same way as generalized Fourier series are introduced.

\section{One dimensional closed model}

Let us consider the standard circle $S^1$ (one dimensional sphere) as our metric space $\Omega$.
This space can be considered as the orbit of a point under the action of the $U(1)$ (or alternatively $SO(2)$) group.
The position
of any person is then given by its polar coordinate $\phi$ ($\phi\in[-\pi,\pi]$) only. The probability $Q$ is the
function of one real variable by means of
\be \label{Qprop1}
Q(\phi_1,\phi_2)=Q(\phi_1-\phi_2)=Q(|\phi_1-\phi_2|) .
\ee
The topology of the problem tells us that $Q$ is a period function with period $2\pi$. This implies we can expand the function $Q$ into
its Fourier series in trigonometric functions in the mentioned interval $[-\pi,\pi]$:
\bdis
Q(\phi)=\sum_{k=-\infty}^\infty a_k e^{ik\phi} .
\edis
Taking into account (\ref{Qprop1}) (and that $Q$ is real-valued) we have that the Fourier  coefficients obey the equation:
$a_{-k}=a_k$. As an example we shall compute the $P(1,b)$ and $P(2,b)$ functions in our case. With respect to
eq. (\ref{p2appint}) and using orthogonality property of trigonometric basis we have after some algebra

\begin{eqnarray} \label{p1ex1}
& &
P(1,b)=
2\pi R(1-Q(b))\sum_{m=-\infty}^\infty a_m^2 \cos(mb) ,
\end{eqnarray}
where the parameter $R$ is introduced as the radius of the considered circle. This parameter will not enter the results because its serves only
as a geometrical tool to keep the distance between any two neighbours equal to $1$ if suitable. This simply means that the volume form is
$R{\rm d}\phi$. Supposing the persons $A$ and $B$ are distanced enough ($b\gg 1 \ \Rightarrow \ Q(b)\ll 1$) we can write instead of
eq. (\ref{p1ex1}) its simpler form:
\be \label{p1ex1a}
P(1,b)\approx 2\pi R\sum_{m=-\infty}^\infty a_m^2 \cos(mb) .
\ee

Using the same idea we can find also the probability $P(2,b)$:
\begin{eqnarray} \label{p2ex1}
& &
P(2,b)=\nonumber \\
& &
(2\pi)^2R^2 (1-Q(b))
\left\{ \sum_{m=-\infty}^\infty a_m^3e^{-imb}-2\sum_{m=-\infty}^\infty\sum_{n=-\infty}^\infty a_m a_n a^2_{m+n}e^{-imb}+ \right. \nonumber \\
& &
\left. \sum_{m=-\infty}^\infty\sum_{n=-\infty}^\infty\sum_{p=-\infty}^\infty a_{m+n+p}a_{m+n}a_m a_n a_p e^{-i(m+p)b} \right\} .
\end{eqnarray}
This expression is a little bit more complicated, assuming, in the same way as previously, $b$ is large enough and neglecting the sums consisting of higher than third power of
the Fourier coefficient we have the approximative expression for the probability $P(2,b)$:
\be \label{p2ex1a}
P(2,b)\approx (2\pi)^2R^2\sum_{m=-\infty}^\infty a_m^3e^{-imb}=(2\pi)^2R^2\sum_{m=-\infty}^\infty a_m^3 \cos(mb) .
\ee
Thus, the important quantity as the clustering coefficient is expressed as follows immediately from eq. (\ref{p2ex1}) putting $b=0$:
\begin{eqnarray} \label{cc1}
& &
\langle C\rangle=\frac{4\pi^2 R^2}{\mcal{N}^2}
\left\{ \sum_{m=-\infty}^\infty a_{m}^3-2\sum_{m=-\infty}^\infty\sum_{n=-\infty}^\infty a_m a_n a^2_{m+n}+ \right. \nonumber \\
& &
\left. \sum_{m=-\infty}^\infty\sum_{n=-\infty}^\infty\sum_{p=-\infty}^\infty a_{m+n+p}a_{m+n}a_m a_n a_p\right\}\approx
\frac{4\pi^2 R^2}{\mcal{N}^2}\sum_{m=-\infty}^\infty a_{m}^3 .
\end{eqnarray}
Furthermore, one could easily generalize the results (\ref{p1ex1a}) and (\ref{p2ex1a}) to the probability $P(k,b)$, the result reads:
\be \label{pkex1}
P(k,b)=(2\pi R)^k\sum_{n=-\infty}^\infty a_n^{k+1}\cos(nb) .
\ee
Now we shall apply our theoretical results to a concrete example of primary probability function $Q$ often discussed in literature in this
context. \\

{\it Uniform distribution within a fixed radius} \\

\begin{figure}
\centering
\includegraphics[width=7cm]{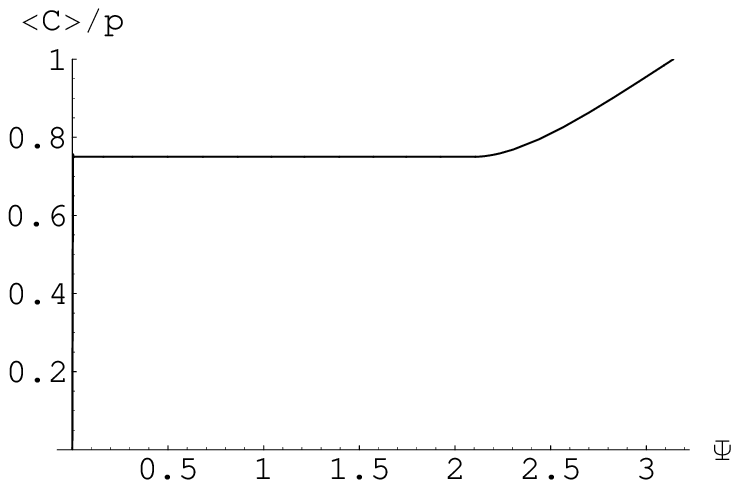}
\includegraphics[width=7cm]{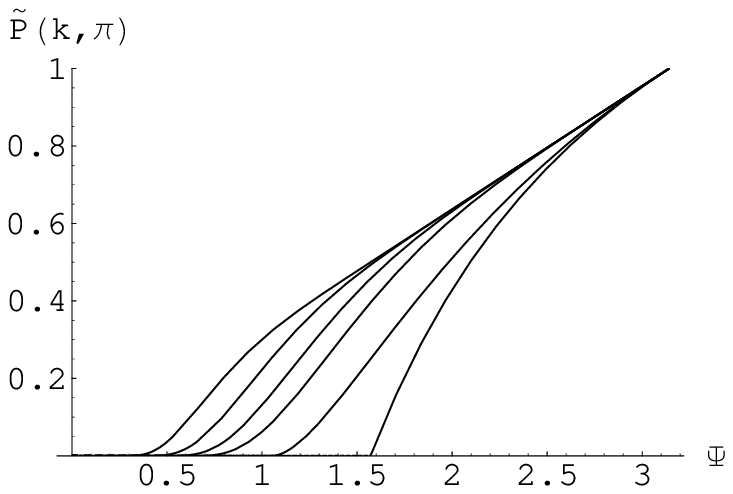}
\caption{Left graph shows the plot of the fraction $\langle C\rangle/p$ for the clustering coefficient given by (\ref{ccexample1}).
The fraction is nearly a constant
(except the value at $\Phi=0$ that is, however, trivial) up to an angle a little bit greater than $2$ and then starts to rise up to $1$ as it
should with respect to (\ref{ccexample1a}). Right one shows the plot of the functions: $\tilde{P}(k,\pi)=P(k,\pi)/(\pi \mcal{N}^k)$
in dependence on the angle $\Phi$. The curves are plotted
(from the right to the left) for $k=1,2,4,6,10,20$. One can obviously identify critical values of the $\Phi$-parameter for $\langle C\rangle/p$ as well as
for $\tilde{P}(k,\pi)$. }
\label{f1}
\end{figure}

Now, the function $Q$ will represent the uniform distribution of probability within a
fixed radius (angle), namely:
\be \label{qexample1}
Q(\phi)=\left\{ \begin{array}{rcl} p &,& \phi\in[-\Phi,\Phi] \\ 0 & , & \phi\notin [-\Phi,\Phi]\end{array} \right. ,
\ee
where $\Phi$ is a parameter lying in the interval $[0,\pi]$ and $p$ is the probability of an acquaintance, therefore $p\in[0,1]$. The mean
number of acquaintances of any given person in this model is given by
\bdis
\mcal{N}=2R\int_0^{\Phi}p{\rm d}\phi=2Rp\Phi .
\edis
In order to obtain the clustering coefficient we are to find the Fourier coefficients of the function (\ref{qexample1}) in the interval
$[-\pi,\pi]$:\ $a_n=\frac{1}{2\pi}\int_{-\pi}^\pi Q(\phi)e^{in\phi}{\rm d}\phi=\frac{1}{\pi}\int_0^\pi Q(\phi)\cos(n\phi){\rm d}\phi$;
\begin{eqnarray*}
& &
a_0=\frac{1}{\pi}\int_0^{\Phi}p{\rm d}\phi=\frac{p\Phi}{\pi} , \quad
a_n=\frac{1}{\pi}\int_0^{\Phi}p\cos(n\phi){\rm d}\phi=\frac{p}{\pi}\frac{\sin(n\Phi)}{n}.
\end{eqnarray*}
With respect to (\ref{cc1}) the clustering coefficient is given by:
\begin{eqnarray} \label{ccexample1}
& &
\langle C\rangle=\frac{p}{\pi\Phi^2}
\left[ \Phi^3+2\sum_{n=1}^\infty\frac{\sin^3(n\Phi)}{n^3}\right] .
\end{eqnarray}
Special case arises when $\Phi=\pi$, i.e. when one considers uniform distribution of probability throughout whole circle. In such a situation we have trivial
result:
\be \label{ccexample1a}
\langle C\rangle(\Phi=0)=p .
\ee
The result (\ref{ccexample1}) is plotted in the graph \ref{f1}. The probabilities $P(k,b)$ in this model are given by (\ref{pkex1}):
\begin{eqnarray}
& &
P(k,b)=(2\pi R)^k\left[\left(\frac{p\Phi}{\pi}\right)^{k+1}+2\sum_{n=1}^\infty \left(\frac{p\sin(n\Phi)}{\pi n}\right)^{k+1}\cos(nb)\right]=
\nonumber \\
& &
\frac{p}{\pi}\left(\frac{\mcal{N}}{\Phi}\right)^k\left[\Phi^{k+1}+2\sum_{n=1}^\infty\frac{\sin^{k+1}(n\Phi)}{n^{k+1}}\cos(nb)\right] ,
\end{eqnarray}
and for the maximally distanced persons:
\be \label{pkpiex1}
P(k,\pi)=
\frac{p}{\pi}\left(\frac{\mcal{N}}{\Phi}\right)^k\left[\Phi^{k+1}+2\sum_{n=1}^\infty\frac{\sin^{k+1}(n\Phi)}{n^{k+1}}(-1)^{n}\right] .
\ee

\section{Higher dimensional models}

Now we shall briefly show how the previous consideration can be extended into higher dimensional case.
For example, people living on Earth are to be
modeled as, at least, two dimensional network (in fact one should incorporate also sociological dimension rather than geographical only). \\

\begin{figure}[h]
\centering
\includegraphics[]{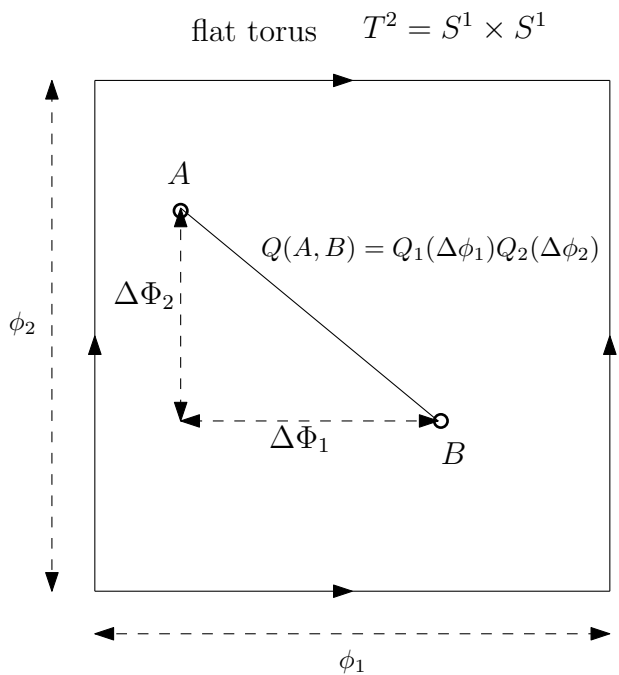}
\caption{Two dimensional flat torus as two dimensional compact background space.}
\label{t2fig}
\end{figure}

There is an straightforward possibility how to generalize our one-dimensional example.
Instead of $S^1$ we can consider the $K$-dimensional tori $T^K$, i.e. the network with the topology
\bdis
T^K:= \underbrace{S^1\times S^1\times \dots \times S^1}_{K-\mbox{times}} .
\edis
The edges between any two neighbours have the unit length. Let $(\phi_1,\phi_2,\dots ,\phi_K)$ be the (periodic) coordinates. The probability
that the persons $A=(0,\dots ,0)$ and $B=(\phi_1,\dots ,\phi_K)$ knows each other is given by the functions of angular distances between the persons
in each direction:
\be \label{qmanyd}
Q(A,B)=Q_1(\phi_1)\cdot Q_2(\phi_2)\cdots Q_K(\phi_K) ,
\ee
where any function $Q_i$ plays the role of single function $Q$ from the previous section. The prescription
(\ref{qmanyd}) for the probability $Q$ allows for introducing an anisotropy in the model, if suitable. This situation in two dimension
is shown in figure
\ref{t2fig}. Let us mention that this example shows we do not need to have metric space, the metric structure is replaced in this case by the
capability to measure distances in separate directions. \\

The task to find the probabilities $P(k,b_1,\dots ,b_K)$ can be obviously solved in the same way as it was done in one
dimensional case in previous section. The difference stands in only one detail that the single integrals
(defining Fourier coefficients) are replaced by multiple ($K$-tuple) integrals defining the Fourier coefficients of
functions defined on $\underbrace{S^1\times S^1\times \dots \times S^1}_{K-\mbox{times}}$. Of course, a practical
difference can appear, namely it is reasonable to expect the numerics in higher dimensional case will be more
complicated. Moreover, one can use the same procedure for any homogenous space with help of generalized Fourier series. Especially, one can
perform analogical computation also in the case when (at least in one dimension) the background space is unbounded - the only thing is to replace
Fourier series by Fourier (generalized) transformation. Generalized orthogonality will allow for analogical expression of the clustering
coefficient also in the case of non-compact symmetry group.

\noi {\bf Acknowledgement:} \ This work was supported partially by GRATEX Research Center and by the VEGA
agency under project no. 1/3042/06 and by the UK grant 403/2007.

\noindent
{\bf Michal Demetrian} is an assistant professor at the department of mathematical and numerical analysis of the Faculty of mathematics, physics and computer science
of the Comenius University in Bratislava. He is mainly interested in critical phenomena in physical dynamical systems interacting with gravity. \\
{\bf Martin Neh\' ez} is an assistant professor at the Faculty of electrical engineering of the Czech Technical University in Prague and PhD. student at the
Comenius University in Bratislava. His research interests are algorithmic graph theory and critical phenomena in discrete structures.


\begin{thebibliography}{25}
\bibitem{milgram}
Milgram, S.: The small world problem, Psychology Today, 1967, 60-67 .
\bibitem{medo}
Medo, M.: Distance-dependent connectivity: Yet Another Approach to the Small World Phenomenon, Physica A, 360/2 (2006), 617-628.
\bibitem{miller}
Miller, G.A. \emph{et al.}: Clustering Coefficients of Protein-Protein Interaction Networks,  Phys. Rev. E 75 (2007), 051910.
\bibitem{lambiotte}
Lambiotte, R.: How does degree heterogeneity affect an order-disorder transition?, EPL, 78 (2007) 68002. Fronczak, P., Fronczak, A. and
Holyst, J.A.: Phase transition in social networks, arXiv: physics/07011182.
\bibitem{reka}
Albert, R. - Barab\' asi, A.L.: Statistical mechanics of complex networks, Rev. Mod. Phys. 74 (2002), 47.
\bibitem{watts}
Watts, D.J.: Small worlds: The Dynamics of Networks between Order and Randomness, Princeton Univ. Press (2003).
\bibitem{watts2}
Watts, D.J. - Strogatz, S.H.: Collective dynamics of 'small-world' networks, Nature 393 (1998), 440-442.
\end{thebibliography}
\end{document}